\documentclass[a4paper,10pt]{article}
\usepackage[utf8]{inputenc}
\usepackage{graphicx}
\usepackage{amsmath,amsfonts,amssymb,graphics,amsthm}

\usepackage{tikz}

\usepackage[colorlinks=false,pdfborderstyle={/W 1}]{hyperref}
\def\arXiv#1#2{\href{http://front.math.ucdavis.edu/#1}{{\tt arXiv:#1 [#2]}}} 
\def\arXivo#1{\href{http://front.math.ucdavis.edu/#1}{{\tt [arXiv:#1]}}} 

\newtheorem{thm}{Theorem}[section]
\newtheorem{quest}[thm]{Question}

\newtheorem{defi}[thm]{Definition}

\title{A unimodular Liouville hyperbolic souvlaki\\ --- an appendix to \arXivo{1603.06712}}
\author{G\'abor Pete \and Gourab Ray}
\date{}

\begin{document}
 \maketitle

\begin{abstract} 
Carmesin, Federici, and Georgakopoulos {\tt [arXiv:1603.06712]} constructed a transient hyperbolic graph that has no transient subtrees and that has the Liouville property for harmonic functions. We modify their construction to get a unimodular random graph with the same properties.
\end{abstract}

\bigskip
\hfill{\small{\it Dedicated to Findus, who has always wanted to grow a tree out of meatballs} \cite{Findus}}

\section{Motivation and the result}

This short note is an appendix to the paper \cite{souvlaki} of Carmesin, Federici, and Georgakopoulos, who constructed a transient hyperbolic graph that has no transient subtrees and that has the Liouville property for harmonic functions. Graphs with such strange properties cannot be transitive. On the one hand, a transitive transient graph has at least 3-dimensional volume growth (as follows from an extension of Gromov's polynomial growth theorem by Trofimov \cite{Trofimov} and Losert \cite{Losert}; see also \cite[Theorem 5.11]{Woess}), hence by the Coulhon-Saloff-Coste isoperimetric inequality (see \cite{CSC} or \cite[Theorem 6.29]{LPbook}) it has at least 3-dimensional isoperimetry, and hence by Thomassen's result \cite{Tho} it contains a transient subtree. On the other hand, a transitive transient hyperbolic graph must be non-amenable \cite{BrHae}, and non-amenable transitive graphs are non-Liouville \cite{KV} and contain non-amenable subtrees \cite{BS:Cheeger}.

This raises the question whether such exceptional graphs can possess any sort of homogeneity. Beyond transitivity, a very natural class of graphs, especially when random walks are considered, is the class of unimodular random graphs, introduced in \cite{BS:limit}, studied in depth in \cite{urn} and in many works since; see \cite[Chapter 14]{PGG} for an overview. Here are the main definitions.

Let $\mathcal{G}_{\star}$ be the space of isomorphism classes of locally finite labeled rooted graphs, and let $\mathcal{G}_{\star\star}$ be the space of isomorphism classes of locally finite labeled graphs with an ordered pair of distinguished vertices, each equipped with the natural local topology: two (doubly) rooted graphs are ``close'' if they agree in ``large'' neighborhoods of the root(s). 

\begin{defi}
We say that a Borel measure $\mu$ on $\mathcal{G}_{\star}$ is \emph{unimodular} if it obeys the \emph{Mass Transport Principle}:
$$ 
\int_{\mathcal{G}_{\star}} \sum_{x\in V(G)}f(G,o,x) \,d\mu(G,o) = \int_{\mathcal{G}_{\star}} \sum_{x\in V(G)} f(G,x,o) \,d\mu(G,o)\,,
$$
for any Borel function $f: \mathcal{G}_{\star\star}\to [0,\infty]$.
\end{defi}

There are several other equivalent definitions; see \cite[Definition 14.1]{PGG}. Probably the nicest one, which works in most situations (e.g., bounded degree non-deterministic graphs), is that the Markov chain on $\mathcal{G}_{\star}$ generated by continuous time random walk on $G$ with rate 1 exponential clocks on the edges is reversible.

An important class of unimodular graphs consists of Cayley graphs of finitely generated groups and of invariant random subgraphs of a Cayley graph. Another one is the class of sofic measures: the closure of the set of finite graphs with a uniform random root under local weak convergence, which is just weak convergence of measures in the space $\mathcal{G}_{\star}$.

Since many results on random walks and harmonic functions on transitive graphs generalize to unimodular or, more generally, stationary random graphs \cite{BC,GLP,AHNR2}, it is natural to ask what the situation is in the present case. Here is our answer:

\begin{thm}
There exists a bounded degree unimodular random graph that is a.s.~transient and hyperbolic, but Liouville and has no transient subtree.
\end{thm}

Our construction will be a splice between the one in the main paper and the so-called $d$-regular canopy tree, which is the local weak limit of larger and larger balls in the $d$-regular tree. It is partly motivated by \cite{Dori}, where similar counterexamples for Bernoulli percolation are constructed based on the canopy tree. However, making the splice is not entirely straightforward here, since we have to put the meatballs on the canopy tree in a way that the graph remains unimodular. For this, the exponentially growing meatballs of the original construction would not work. 

Having seen the first version of this paper, Itai Benjamini asked the following:

\begin{quest}
Is there a bounded degree unimodular random graph that is non-Liouville but contains no transient subtree?
\end{quest}

A non-unimodular example is the integer line with a copy of a transient Souvlaki attached to each vertex. Another, planar, example is the binary tree, with each edge at level $n$ replaced by two binary trees of depth $3^n$ meeting at the leaves like this: $< >$. However, we have been unable to make unimodular versions of these examples or to prove that they cannot exist. On the other hand, our understanding is that Tom Hutchcroft has recently managed to construct a unimodular example.

\section{The construction}

Take a $d$-ary tree $T_n$ of height $n$ (i.e., the root has $d$ children, each of which has $d$ children, and so on, stopping with the $n$th descendent generation). We are going to replace each
edge of $T_n$ by a modification of the graphs $M_k$ in the Souvlaki
construction of the main paper. Recall what $W,H_2,H_3$ are.
Consider a subpath of the bottom double ray $P_k$ of $W$ of length
$(k-1)^2+k^2+k^4$. Define $M_k$ to be the subgraph of $H_3$ induced by vertices of the form $(t,w)$
such that $w$ lies at or above $P_k$ and has height $h(w)$ at most
$k$. We call $M_k$ the meatballs.

We now ``replace'' each edge of $T_n$ at height $n-k+1$ (edges such that the vertex closer to the root has height $n-k$) by $M_k$ for $k=1,\ldots,n$. The word replace is within quotes because we have to specify the way we glue adjacent meatballs. We
divide $P_k = L_k  \cup R_k \cup A_k$, where $L_k$ is the segment of the leftmost $k^2$ vertices, $R_k$ is the segment
of the rightmost  $(k-1)^2$ vertices, and $A_k$ is the middle $k^4$ vertices. Now let
$M^L_k$ denote the set of vertices $(t,w)$ so that $w$ lies on or
above $L_k \cup A_k$ and has height
at most  $k$. Define $M^R_k$ to be the graph induced by the rest of the vertices. Note that $M^R_k$ and $M^L_k$
are joined together by a set of edges. Let $B_k$ denote the endpoint of these edges that lie in $M_k^L$. 

Since we have the tree $T_n$ instead of just a line, we need to modify $M_k$ a bit, so that it branches into $d$ copies of $M^R_k$ for the identifications.
That is, we take one copy of $M_k^L$ and $d$ copies of $M^R_k$, then glue each of the latter with $M_k^L$ along $B_k$. Call this new gadget $M_k'$. The height function $h$ extends to $M_k'$, with values between $0$ and $k$.

Now take an edge $e$ at height $n-k+1$, for $k \ge 1$. Remove it and replace it with $M_{k}'$ so that $e$ corresponds to the segment
$L_k \cup A_k$, while its $d$ children $e_1,\ldots e_d$ at height $n-k+2$ correspond to the $d$ copies of $R_k$ for $e$. Since each $R_k$ contains $(k-1)^2$ vertices, we can identify them with the copies of $L_{k-1}$ for the edges $e_i$. This completes the gluing procedure. For an edge $e$ at height
$n$, replace $e$ just by $M^L_1$, without branching into copies of $M^R_1$. Call the new graph so obtained $T_n'$.

Now pick a uniform random vertex $\rho_n$ from $T_n'$ and take a weak limit. Call the limit $(T,\rho)$. Clearly this graph is unimodular,
Gromov hyperbolic and bounded degree, from arguments in the main paper.

\section{Proofs}

\paragraph{Root height.} Pick any integer $d>6$. At height $i$ of $M^L_{k}$, the number of vertices is $3^i 2^i (k^2+k^4)$, with the factor $3^i$ coming from  $H_2$, and the factor $2^i$ coming from $W$. Thus, the volume of $M^L_{k}$ is 
$$v_k:=\frac{6^{k+1}-1}{5} (k^4+k^2),$$ 
and the probability that the uniform root in $T_n$ is a vertex in one of the $M^L_{k}$s is
$$
p_{k,n} := \frac{v_k d^{n-k+1}}{\sum_{j=1}^n v_jd^{n-j+1}} = \frac{v_k d^{-k}}{\sum_{j=1}^n v_j d^{-j}}\,.
$$
Since $d>6$ implies that $v_kd^{-k}$ is summable, the limit $p_k:=\lim_{n\to\infty} p_{k,n}$ is a proper probability distribution for $k=1,2,\dots$. Therefore, the root in the limit $(T,\rho)$ is almost surely at a level corresponding to a finite $k$, at a finite distance from the leaves (of the underlying canopy tree that is the local weak limit of the original trees $T_n$). In other words, we can think of the limit $(T,\rho)$ as a souvlaki with a canopy tree skewer, with a random root somewhere.


\paragraph{Constructing a good flow.} One can think of the canopy tree as an infinite spine with finite bushes hanging off of it. Similarly, our canopy tree souvlaki has an infinite spine, a ``traditional'' infinite souvlaki. It is of course enough to show that this infinite spine is transient. We will construct for each $k\geq 1$ a unit flow $g$ from $R_{k+1}$ to $L_{k+1}$, with an energy that is summable in $k$. Concatenating these flows yields a flow along the spine to infinity, with finite energy, hence the spine turns out to be transient by Terry Lyons' criterion. (Unfortunately, the roles of $R$s and $L$s are now swapped compared to the main paper, due to the way that the infinite limit is constructed.) 

Note that $R_{k+1}$ has $k^2$ vertices and $L_{k+1}$ has $(k+1)^2$ vertices. We name the vertices
in $R_{k+1}$ as  $r_1,r_2,\ldots, r_{k^2}$ and the vertices in $L_{k+1}$ as $l_1, \ldots, l_{(k+1)^2}$. 
We will construct a flow $g_j$ from $r_j$
to $L_{k+1}$ and $g$ will be the sum the flows $g_j$. The flow $g$
will have outflow $1/k^2$ from each vertex in $R_{k+1}$ and inflow
$1/(k+1)^2$ for each vertex in $L_{k+1}$.

In the main paper, since the meatballs had exponentially growing lengths,
there was a straightforward division of the total flow from each
vertex into two vertices. Since here everything is polynomial, the division
is slightly more complicated, but this is just a technicality. We deal with this
as follows.

Since $H_2$ is transient, there exists a natural unit flow $t$ (equally branching off at each vertex) with finite energy. Then, 
the flow $g_j$ is just $\frac 1 {k^2} t$ up to height $k$, for each $j=1,\ldots, k^2$. After this, $g_j$ takes $k^2/(k+1)^2$ fraction of the total
incoming flow at height $k$, flows along the horizontal edges at height $k$ until reaching above $l_j$, then flows along the tree proportionally with $t$
(in reverse direction) to reach $l_j$. So, for $j=1,\ldots, k^2$, the total flow into $l_j$ is already $1/(k+1)^2$.

There is still $(1-k^2/(k+1)^2)$ fraction of the outflow from $R_{k+1}$ left at level $k$. In fact, there is $(k^{-2}-(k+1)^{-2} )3^{-k}$ amount left in
each vertex at height $k$, in each $g_j$, $1 \le j \le k^2$. Flow this amount horizontally along level $k$ until we reach above the vertex $l_{k^2+1}$. The total flow here (summed over $j$) is now $((k+1)^2-k^2)/(k+1)^2$. Now we take $1/(k+1)^2$ out of this, and flow it along the tree above $l_{k^2+1}$ proportionally with $t$, in reverse direction. We flow the remaining amount of flow at level $k$ horizontally to reach above $l_{k^2+2}$, and again drop $1/(k+1)^2$ amount along the tree with $t$. We continue like this until all the flow is exhausted. Note that we input a flow $1/(k+1)^2$ for each $l_j$, $k^2< j\le (k+1)^2$. Thus all in all, this is a unit flow from $R_{k+1}$ to $L_{k+1}$.

\paragraph{Energy computation.}
The total energy cost for going up the tree to height $k$ for each $r_j$ is $E(t)/k^4$. The flow received at
each vertex is $3^{-k}/k^2$. On each horizontal edge at height $k$, there are at most $k^2$ flows that we are summing up, hence 
total flow through is at most $3^{-k}$. The total length of horizontal
paths is $O(k^4) 3^k 2^k$, hence the total energy along the horizontal edges is $O(k^4) 6^k 3^{-2k}$.
The energy for going down each tree above $l_j$ is again at most $E(t)/k^4$. Therefore, the total energy dissipation is 
$$
O(k^4) 6^k 3^{-2k} + O(k^2) E(t)/k^4 = O(k^4) (2/3)^k + O(1/k^2) = O(1/k^2)\,,
$$
which is summable in $k$. This concludes the proof of transience.


\paragraph{Liouville property.} Removing the infinite spine, the canopy tree souvlaki falls apart into finite pieces. Thus, random walk started anywhere in the graph will almost surely hit the spine. This and the Optional Stopping Theorem for bounded martingales imply that any bounded harmonic function is determined by its restriction to the spine. Also, the radial symmetry of the harmonic function in the meatballs along the spine is preserved, as in the main paper. Thus the graph is Liouville by the same argument as in the main paper.


\paragraph{Transient subtree.} We first claim that any subtree inside the infinite spine must be recurrent. This can be shown by the same argument as in the main paper. Namely, in the proof of Theorem 6.1, the degrees $d(v_k)$ and $d(v_{k+1})$ are $O(k^2)$, hence following the proof of Lemma 6.2 shows that taking $s=k^4$ is enough for the effective resistance between $v_k$ and $v_{k+1}$ to be uniformly positive. This is exactly the choice we made in defining $A_k$, thus the effective resistance of any subtree to infinity is infinite, hence it is recurrent.
	
	Now, for a general subtree, when we do the contraction into the vertices $v_k$, then the portions of the subtree inside the finite bushes off the spine get contracted into finite pieces, each attached to the contracted graph from the spine at a single vertex. These finite pieces do not influence transience of the contracted graph, hence the original subtree is also recurrent.

\paragraph{Acknowledgments.} We are grateful to Itai Benjamini and Agelos Georgakopoulos for comments on the manuscript. 


\vskip 0.3 cm

\noindent {\bf G\'abor Pete}\\
R\'enyi Institute, Hungarian Academy of Sciences, Budapest, and\\
Institute of Mathematics, Budapest University of Technology and Economics\\
\url{http://www.math.bme.hu/~gabor}\\
Partially supported by the Hungarian National Research, Development and Innovation Office, NKFIH grant K109684, and by the MTA R\'enyi Institute ``Lend\"ulet'' Limits of Structures Research Group.
\vskip 0.3 cm

\noindent {\bf Gourab Ray}\\
Statistical Laboratory, University of Cambridge\\
\url{https://sites.google.com/site/gourabmathematics/}
The work supported by the Engineering and Physical Sciences Research Council under grant EP/103372X/1.

\end{document}